\documentstyle[12pt, twoside]{article}
\input amssym.def
\input amssym.tex

\newenvironment{demo}[1]%
{\vskip-\lastskip\medskip
  \noindent
  {\em #1.}\enspace
  }%
{\qed\par\medskip
  }

\newcommand{\qed}{
  \strut\hfill
  \mbox{$\Box$}
  }

\newtheorem{theorem}{Theorem}
\newtheorem{lemma}{Lemma}
\newtheorem{remark}{Remark}

\newtheorem{proposition}{Proposition}
\newtheorem{corollary}{Corollary}

\newcommand{\az}{a(z)}
\newcommand{\C}{\Bbb C}
\newcommand{\CG}{\C [\G]}
\newcommand{\ch}{\mbox{ch}}
\newcommand{\De}{\Delta}
\newcommand{\Dec}{\Delta_c}

\newcommand{\diffi}{\frac{\partial}{\partial p_i} }
\newcommand{\diffij}{\frac{\partial}{\partial p_{i +j}} }
\newcommand{\diffj}{\frac{\partial}{\partial p_j} }

\newcommand{\diffig}{\frac{\partial}{\partial p_i (\g)} }
\newcommand{\diffijg}{\frac{\partial}{\partial p_{i +j} (\g)} }
\newcommand{\diffjg}{\frac{\partial}{\partial p_j (\g)} }

\newcommand{\f}{f^{\lambda}}
\newcommand{\fbe}{f^{\lambda(\beta)}}
\newcommand{\fge}{f^{\lambda(\gamma)}}

\newcommand{\fh}{\frak{h}}
\newcommand{\G}{\Gamma}
\newcommand{\Gm}{{\Gamma}_m}
\newcommand{\Gn}{{\Gamma}_n}
\newcommand{\g}{\gamma}
\newcommand{\hf}{  \frac12  }
\newcommand{\hg}{\widehat{\frak h}_{\G}}         
\newcommand{\hh}{\widehat{\frak h}}
\newcommand{\la}{\lambda}
\newcommand{\laz}{\lambda_c}
\newcommand{\Lz}{L(z)}
\newcommand{\be}{\beta}

\newcommand{\RG}{R_{ \G}}
\newcommand{\Rz}{R_{\Bbb Z} (\G)}
\newcommand{\schu}{s_{\lambda}}
\newcommand{\sbe}{s_{\lambda(\beta)} (\beta)}
\newcommand{\sge}{s_{\lambda(\gamma)} (\gamma)}
\newcommand{\SL}{S_L}
\newcommand{\SG}{ S_\G }

\newcommand{\Z}{{\Bbb Z}}
\begin{document}

\title{Virasoro algebra and wreath product convolution}

\author{Igor B.~Frenkel and Weiqiang Wang
}
\date{}


\maketitle
\begin{abstract}
We present a group theoretic construction of the Virasoro algebra
in the framework of wreath products. This can be regarded as a
counterpart of a geometric construction of Lehn in the theory of
Hilbert schemes of points on a surface.
\end{abstract}

\section*{Introduction}
It is by now well known that a direct sum $\bigoplus_{n \geq 0}
R(S_n)$ of the Grothendieck rings of symmetric groups $S_n$ can be
identified with the Fock space of the Heisenberg algebra of rank
one. One can construct vertex operators whose components generate
an infinite-dimensional Clifford algebra, the relation known as
boson-fermion correspondence \cite{F} (also see \cite{J}). A
natural open problem which arises here is to understand the group
theoretic meaning of more general vertex operators in a vertex
algebra \cite{B, FLM}.

A connection between a direct sum $\RG = \bigoplus_{n \geq 0}
R(\G_n)$ of the Grothendieck rings of wreath products $\Gn = \G
\sim S_n$ associated to a finite group $\G$ and vertex operators
has been realized recently in \cite{W} and \cite{FJW} (also see
\cite{Z, M} for closely related algebraic structures on $\RG$).
When $\G$ is trivial one recovers the above symmetric group
picture. On the other hand, this wreath product approach turns out
to be very much parallel to the development in the theory of
Hilbert schemes of points on a surface (cf. \cite{W, N2} and
references therein). As one expects that new insight in one theory
will shed new light on the other, this refreshes our hope of
understanding the group theoretic meaning of general vertex
operators.

The goal of this paper is to take the next step in this direction
to produce the Virasoro algebra within the framework of wreath
products. Denote by $\G_*$ the set of conjugacy classes of $\G$,
and by $c^0$ the identity conjugacy class. Recall \cite{M, Z} that
the conjugacy classes of the wreath product $\Gn$ are
parameterized by the partition-valued functions on $\G_*$ of
length $n$ (also see Section~\ref{sect_wreath}). Given $c \in
\G_*$, we denote by $\laz $ the function which maps $c$ to the
one-part partition $(2)$, $c^0$ to the partition $( 1^{n-2})$ and
other conjugacy classes to $0$. We will define an operator $\Dec$
in terms of the convolution with the characteristic class function
on $\Gn$ (for all $n$) associated to the conjugacy class
parameterized by $\laz$. We show that this operator can be
identified with a differential operator which is the zero mode of
a certain vertex operator when we identify $\RG$ as in \cite{M}
with a symmetric algebra. A group theoretic construction of
Heisenberg algebra has been given in \cite{W} (also see
\cite{FJW}) which acts on $\RG$ irreducibly. The commutator
between $\Dec$ and the Heisenberg algebra generators on $\RG$
provides us the Virasoro algebra generators.

Our construction is motivated in part by the work of Lehn \cite{L}
in the theory of Hilbert schemes. Among other results, he showed
that an operator defined in terms of intersection with the
boundary of Hilbert schemes may be used to produce the Virasoro
algebra when combined with earlier construction of Heisenberg
algebra due to Nakajima and independently Grojnowski \cite{N1,
Gr}. It remains an important open problem to establish a precise
relationship between Lehn's construction and ours.

We remark that the convolution operator in the symmetric group
case (i.e. $\G$ trivial), when interpreted as an operator on the
space of symmetric functions, is intimately related to the
Hamiltonian in Calegero-Sutherland integrable system  and to the
Macdonald operator which is used to define Macdonald polynomials
\cite{AMOS, M}.

After we discovered our group theoretic approach toward the
Virasoro algebra, we notice that our convolution operator in the
symmetric group case has been considered by Goulden \cite{Go} when
studying the number of ways of writing permutations in a given
conjugacy class as products of transpositions. We regard this as a
confirmation of our belief that the connections between $\RG$ and
(general) vertex operators are profound. It is likely that often
time when we understand something deeper in this direction, we may
realize that it is already hidden in the vast literature of
combinatorics, particularly on symmetric groups and symmetric
functions for totally different needs. Then the virtue of our
point of view will be to serve as a unifying principle which
patches together many pieces of mathematics which were not
suspected to be related at all.

The plan of this paper is as follows. In Sect.~\ref{sect_vertex},
we recall some basics of Heisenberg and Virasoro algebras from the
viewpoint of vertex algebras. In Sect.~\ref{sect_wreath} we set up
the background in the theory of wreath products which our main
constructions are based on. In Sect.~\ref{sect_conv} we present
our main results. Some materials in the paper are fairly standard
to experts, but we have decided to include them in hope that it
might be helpful to the reader with different backgrounds.
\section{Basics of Heisenberg and Virasoro algebras} \label{sect_vertex}
In the following we will present some basic constructions in
vertex algebras which give us the Virasoro algebra from a
Heisenberg algebra.

Let $L$ be a rank $N$ lattice endowed with an integral
non-degenerate symmetric bilinear form. Indeed we will only need
the case $L = \Z^N$ with the standard bilinear form.

Denote by $\fh = \C \bigotimes_{\Z} L$ the vector
space generated by $L$ with the bilinear from $\langle-,-\rangle$
induced from $L$. We define the Heisenberg algebra
$$
  \hh = \C [t, t^{-1}] \bigotimes \fh \bigoplus \C C
$$
with the following commutation relations:
\begin{eqnarray*}
   [ C, a_n ]     & =& 0  , \\
   {[ a_n, b_m ]} & =& n \delta_{n, -m} \langle a,b \rangle C,
\end{eqnarray*}
where $a_n $ denotes $t^n \otimes a, a \in \fh$.

We denote by $ \SL$ the symmetric algebra generated by
$\hh^{-} = t^{-1} \C [t^{-1}] \bigotimes \fh$.
It is well known that $\SL$ can be given the structure
of an irreducible module over Heisenberg algebra $\hh$
by letting $a_{-n}, n>0$ acts as multiplication
and letting
\begin{eqnarray*}
  a_n . a_{-n_1}^1 a_{-n_2}^2 \ldots a_{-n_k}^k
  = \sum_{i =1}^k \delta_{n,n_i} \langle a, a^i \rangle
      a_{-n_1}^1 a_{-n_2}^2 \ldots
      \check{a}_{-n_i}^i \ldots a_{-n_k}^k ,
\end{eqnarray*}
where $n \geq 0, n_i > 0, a, a^i \in \fh$ for $i =1, \ldots , k$,
and $\check{a}_{-n_i}^i $ means the very term is deleted.
A natural gradation on $\SL$ is defined by letting
$$
 \deg (a_{-n_1}^1 a_{-n_2}^2 \ldots a_{-n_k}^k)
  = n_1+  \ldots  +n_k.
$$ We say an operator on $\SL$ is of degree $p$ if it maps any
$n$-th graded subspace of $\SL$ to $(n +p)$-th graded subspace.

The space $\SL$ carries a natural structure of a vertex algebra
\cite{B, FLM}. It is convenient to use the generating function in
a variable $z$: $$
 \az = \sum_{n \in \Z} a_n z^{ -n -1}, \quad a \in \fh.
$$
In the language of vertex algebras, this is the vertex operator
associated to the vector $a_{-1} \in \SL$.
The normally ordered product between two vertex operators
$a(z)$ and $b(z)$ is defined as
\[
  : a(z) b(z) : =
  \sum_{n <0} a_n z^{ -n -1} b(z)
  + b(z) \sum_{n \geq 0} a_n z^{ -n -1}.
\]
We remark that the commutation relations between $a_n, b_m, n, m
\in \Z$ are encoded in the following operator product expansion
(OPE) (cf. \cite{FLM}):
\[
 a(z) b(w) \sim \frac{ \langle a, b \rangle }{ (z -w)^2}.
\]

We recall that the Virasoro algebra is spanned by
$L_n, n \in \Z$ and a central element $C$, with
the following commutation relations:
\[
  [ L_n, L_m] = (n -m) L_{n +m} + \frac{n^3 -n}{12} \delta_{n,m} C.
\]
It is convenient to denote
\[
 \Lz = \sum_{n \in \Z} L_n z^{ -n -2}.
\]
Given an orthonormal basis $ a^i, i =1, \ldots, n$ of $\fh$, we
define a series of operators $\widetilde{L}^i_n$ $ (i =1, \ldots
N, n \in \Z) $ acting on $\SL$ by
\begin{eqnarray*}
 \widetilde{L}^i (z)
     \equiv \sum_{n \in \Z} \widetilde{L}^i_n z^{ -n -2}
                   =  \hf : a^i (z) a^i (z)  : .
\end{eqnarray*}

The following proposition is well known (cf. e.g. \cite{FLM}).

\begin{proposition}
The operators $\widetilde{L}^i_n$ $(i =1, \ldots N, n \in \Z)$
generate $N$ commutative copies of Virasoro algebra of central
charge $1$. The operators $\widetilde{L}_n, n \in \Z $ generate
the Virasoro algebra with central charge $N$. Namely we have
  \begin{eqnarray*}
   {[\widetilde{L}^i_n, \widetilde{L}^j_m  ]}
     = \delta_{ij} (n -m) \widetilde{L}^i_{n +m}
     + \frac{n^3 -n}{12} \delta_{ij} \delta_{n,m} .
 \end{eqnarray*}
\end{proposition}
From now on we will simply write $\widetilde{L}^i_n$ as $L^i_n$.

We introduce operators
$$
 \De^i = \frac16 \int :a^i (z)^3: z^2 dz.
$$
These are well-defined
operators of degree $0$ acting on $\SL$. One can
easily check that
$$
 \De^i =\frac12 \sum_{n, m > 0}
  ( a^i_{ -n -m} a^i_n a^i_m + a^i_{-n}a^i_{-m}a^i_{n+m} ).
$$ Here we omit on the right hand side the terms involving $a_0$
since $a_0$ acts as $0$ on the Fock space $\SL$.

\begin{lemma}  \label{lem_ham}
   We have
 \[
  [\De^i, a_n ] = -n \langle a^i, a \rangle  L^i_n.
 \]
\end{lemma}

\begin{demo}{proof}
 It is clear that $[\De^i, a^j_n ] = 0$ for $i \neq j$.
 So it suffices to prove that
 \[
  [\De^i, a^i_n ] = -n L^i_n.
 \]

     %
     %
We calculate that
 \begin{eqnarray*}
   {[\De^i, a^i (z) ]} & =&
    {\rm Res}_{w=0} \frac16 [ :a^i ( w)^3:, a^i (z) ] w^2  \\
   & =& - \hf {\rm Res}_{w=0}:a^i ( w)^2: w^2 \sum_{m \in \Z}  m w^{m -1} z^{-m -1} \\
   & =& -{\rm Res}_{w=0} \sum_{n \in \Z} L^i _n w^{ -n -2} w^2
        \sum_{m \in \Z}  m w^{m -1} z^{-m -1} \\
   & =& -\sum_{n \in \Z} n L^i _n z^{-n -1},
 \end{eqnarray*}
 Therefore
 $$
  [\De^i, a^i_n ] = {\rm Res}_{z=0} [\De^i, a^i (z)] z^n  = -n L^i _n.
 $$
\end{demo}
\section{Representation rings of wreath products}  \label{sect_wreath}
%
%
%
%
%
%
  Given a finite group $\G$, we denote by $\G^*$ the set of complex
irreducible characters and by $\G_*$ the set of conjugacy classes.
We denote by $\Rz$ the $\Z$-span of irreducible characters of $
\G$. Denote by $c^0$ the identity conjugacy class. We identify
$R(\G) = \C \otimes_{\Z } R_{\Z}(\G)$ with the space of class
functions on $\G$.

For $c \in \G_*$ let $\zeta_c$ be the order of the centralizer of
an element in the class $c$. Denote by $|\G |$ the order of $\G$.
The usual bilinear form $\langle -, - \rangle_{\G}$ on $R(\G )$ is
defined as follows (often we will omit the subscript $\G$):
\begin{eqnarray}  \label{eq_bilin}
 \langle f, g \rangle = \langle f, g \rangle_{\G}
    = \frac1{ | \G |}\sum_{x \in \Gamma}
          f(x) g(x^{ -1})
 = \sum_{c \in \Gamma_*} \zeta_c^{ -1} f(c) g(c^{ -1}),
\end{eqnarray}
where $c^{ -1}$ denotes the conjugacy class
$\{ x^{ -1}, x \in c \}$. Clearly $\zeta_c = \zeta_{c^{-1}}$.
It is well known that
\begin{eqnarray}
  \langle \g, \g' \rangle
   &= & \delta_{\g, \g'}, \quad \g, \g' \in \G^*\nonumber \\
  \sum_{ \g \in \G^*} \g (c ')  \g ( c^{ -1})
    &= & \delta_{c, c '} \zeta_c, \quad c, c ' \in \G_*. \label{eq_stand}
\end{eqnarray}

One may regard $\CG$ as the space of functions on $\G$, and thus
$R(\G)$ as a subspace of $\CG$. Given $f, g \in \CG$, the  {\em
convolution} $f * g \in\CG$ is defined by
\[
  f * g (x) =
    \sum_{y \in \G} f(x y^{-1}) g(y),
  \quad f, g \in \CG, x \in \G.
\]
In particular if $f, g \in R(\G)$, then so is $f *g$. It is well
known that
\begin{eqnarray}  \label{eq_idem}
   \g ' * \g = \frac{|\G|}{d_{\g}} \delta_{\g ,\g '} \g,
    \quad \g', \g \in \G^*,
\end{eqnarray}
where $d_{\g}$ is the {\em degree} of the irreducible character
$\g$.

Denote by $K_c$ the sum of all elements in the conjugacy class
$c$. By abuse of notation, we also regard $K_c$ the class function
on $\G$ which takes value $1$ on elements in the conjugacy class
$c$ and $0$ elsewhere. It is clear that $K_c, c \in \G_*,$ form a
basis of $R(\G)$. The elements $K_c, c \in \G_*$ actually form a
linear basis of the center in the group algebra $\CG$. But we will
not need this fact.
%
%
%
%
%
%
%

For wreath products we basically follow the excellent presentation
of Macdonald \cite{M}, Appendix B, Chapter 1, with the exception
of Theorem \ref{th_heis} which is quoted from \cite{W} (also see
\cite{FJW}). Given a positive integer $n$, let $\Gamma^n = \Gamma
\times \cdots \times \Gamma$ be the $n$-th direct product of
$\Gamma$. The symmetric group $S_n$ acts on $\Gamma^n$ by
permutations: $\sigma (g_1, \cdots, g_n)
  = (g_{\sigma^{ -1} (1)}, \cdots, g_{\sigma^{ -1} (n)}).
$
The wreath product of $\Gamma$ with $S_n$ is defined to be
the semi-direct product
$$
 \Gamma_n = \{(g, \sigma) | g=(g_1, \cdots, g_n)\in {\Gamma}^n,
\sigma\in S_n \}
$$
 with the multiplication
$$
(g, \sigma)\cdot (h, \tau)=(g \, {\sigma} (h), \sigma \tau ) .
$$

Let $\la=(\la_1, \la_2, \cdots, \la_l)$ be a partition of integer
$|\la|=\la_1+\cdots+\la_l$, where $\la_1 \geq \dots \geq \la_l
\geq 1$. The integer $l$ is called the {\em length} of the
partition $\la $ and is denoted by $l (\la )$. We will identify
the partition $(\la_1, \la_2, \cdots, \la_l)$ with $(\la_1, \la_2,
\cdots, \la_l, 0, \cdots, 0)$. We will also make use of another
notation for partitions: $$ \la=(1^{m_1}2^{m_2}\cdots) , $$ where
$m_i$ is the number of parts in $\la$ equal to $i$.

For a finite set $X$ and $\rho=(\rho(x))_{x \in X}$ a family of
partitions indexed by $X$, we write $$\|\rho\|=\sum_{x  \in
X}|\rho(x)|.$$ Sometimes it is convenient to regard
$\rho=(\rho(x))_{x \in X}$ as a partition-valued function on $X$.
We denote by ${\cal P}(X)$ the set of all partitions indexed by
$X$ and by ${\cal P}_n (X)$ the set of all partitions in ${\cal P}
(X)$ such that $\|\rho\| =n$.

The conjugacy classes of ${\Gamma}_n$ can be described in the
following way. Let $x=(g, \sigma )\in {\Gamma}_n$, where $g=(g_1,
\ldots , g_n) \in {\Gamma}^n,$ $ \sigma \in S_n$. The permutation
$\sigma $ is written as a product of disjoint cycles. For each
such cycle $y=(i_1 i_2 \cdots i_k)$ the element $g_{i_k} g_{i_{k
-1}} \cdots g_{i_1} \in \Gamma$ is determined up to conjugacy in
$\Gamma$ by $g$ and $y$, and will be called the {\em
cycle-product} of $x$ corresponding to the cycle $y$. For any
conjugacy class $c$ and each integer $i \geq 1$, the number of
$i$-cycles in $\sigma$ whose cycle-product lies in $c$ will be
denoted by $m_i(c)$. Denote by $\rho (c)$ the partition $(1^{m_1
(c)} 2^{m_2 (c)} \ldots )$, $c \in \G_*$. Then each element $x=(g,
\sigma)\in {\Gamma}_n$ gives rise to a partition-valued function
$( \rho (c))_{c \in \G_*} \in {\cal P} ( \G_*)$ such that
$\sum_{i, c} i m_i(c) =n$. The partition-valued function $\rho =(
\rho(c))_{ c \in G_*} $ is called the {\em type} of $x$. It is
known that any two elements of ${\Gamma}_n$ are conjugate in
${\Gamma}_n$ if and only if they have the same type.

Given a partition $\lambda = (1^{m_1} 2^{m_2} \ldots )$,
we define
$
  z_{\la } = \prod_{i \geq 1}i^{m_i}m_i!.
$
We note that $z_{\la }$ is the order of the centralizer of an
element of cycle-type $\la $ in $S_{|\la |}$. The order of the
centralizer of an element $x = (g, \sigma) \in {\Gamma}_n$ of the
type $\rho=( \rho(c))_{ c \in \G_*}$ is
\begin{eqnarray}   \label{eq_centr}
Z_{\rho}=\prod_{c\in \G_*}z_{\rho(c)}\zeta_c^{l(\rho(c))}.
\end{eqnarray}
%
%
%
%
%
%
Recall that $R(\Gn) = R_{\Z} (\Gn) \otimes_{\Bbb Z} \Bbb C$. We
set
\begin{eqnarray*}
  \RG = \bigoplus_{n \geq 0} R (\Gn).
\end{eqnarray*}

A symmetric bilinear form on $\RG$ is given by
\begin{eqnarray}   \label{eq_series}
\langle u, v \rangle
 = \sum_{ n \geq 0} \langle u_n, v_n \rangle_{\Gn } ,
\end{eqnarray}
where
$u = \sum_n u_n$ and $v = \sum_n v_n$ with $u_n, v_n\in \Gn$.

Since $\Rz$ may be regarded as an integral lattice with
non-degenerate symmetric bilinear form given by (\ref{eq_bilin})
and with an orthonormal basis $\G^*$, we can apply the
constructions in Sect.~\ref{sect_vertex} to the lattice $\Rz$. We
denote the corresponding Heisenberg algebra by $\hg$ with
generators $a_n (\g ), n \in \Z, \g \in \G^*$, and its irreducible
representation by $\SG$.

By identifying $a_{-n} (\g)$ $(n >0, \g \in \G^*)$ with the $n$-th
power sum in a sequence of variables parameterized by $\g \in
\G^*$, we may regard $\SG$ as the algebra of symmetric functions
parameterized by $\g \in \G^*$. In particular the operator $a_n
(\g)$ $(n >0, \g \in \G^*)$ acts as the differential operator $n
\frac{\partial}{\partial a_{-n}(\g)}$.

For $m \in \Bbb Z, c \in \G_*$ we define
\begin{eqnarray}  \label{eq_change}
 a_{ m}( c) = \sum_{ \g\in \G^*} \gamma(c^{-1}) a_m( \g ).
\end{eqnarray}
From the orthogonality of the irreducible characters
(\ref{eq_stand}) it follows that
\begin{eqnarray*}
  a_m( \g )
   = \sum_{c \in \G_*} \zeta_c^{ -1}
   \g  (c) a_m(c).
\end{eqnarray*}
Thus $ a_{n}(c)$ ($ n\in \Bbb Z, c \in \G_*$) and $C$ form a new
basis for the Heisenberg algebra $\hg$.

Given a partition $ \la = ( \la_1, \la_2, \ldots )$
and $c \in \G_*$, we define
\begin{eqnarray*}
  p_{ \la } (c )  = p_{ \la_1}(c) p_{ \la_2 } (c) \cdots.
\end{eqnarray*}
For any $\rho = ( \rho (c) )_{ c \in \G_* } \in
 {\cal P } ( \G_* )$, we further define
\begin{eqnarray*}
 P_{ \rho} = \prod_{ c \in \G_*} p_{ \rho (c)} (c).
\end{eqnarray*}
The elements $P_{ \rho} ,
\rho \in {\cal P}(\G_* )$ form a $ \C$-basis for $\SG $.

We define a bilinear form
$\langle \ ,  \ \rangle $ on the space $ \SG $ by letting
\begin{eqnarray}   \label{eq_sg}
  \langle P_{\rho}, P_{\sigma} \rangle
  = \delta_{\rho, \sigma} Z_{\rho}.
\end{eqnarray}

Let $\Psi : \Gn \rightarrow \SG$ be the map defined by $\Psi (x) =
P_{ \rho}$ if $x \in \Gn$ is of type $\rho$. We define a $\Bbb
C$-linear map $\ch: \bigoplus_{ n \geq 0} \C [\Gn] \longrightarrow
\SG$ by letting
\begin{eqnarray}
 \ch (f )
 & =& \frac 1{| \Gn|} \sum_{x \in \Gn} f(x) \Psi (x)
                              \label{eq_charmap}   \\
 & =& \sum_{\rho \in {\cal P} (\G_*)}
  Z_{\rho}^{-1} f_{\rho} P_{  \rho } \quad \mbox{if } f \in \RG,
    \nonumber
\end{eqnarray}
where $f_{\rho}$ is the value of $f$ at elements of type $\rho$.
Most often we will think of $\ch$ as a map from $\RG$ to $\SG$ as
in \cite{M}, usually referred to as the {\em characteristic map}.
It is well known that $\ch : \RG \rightarrow \SG$ is an isometry
for the bilinear forms on $\RG$ and $\SG$ defined in
(\ref{eq_series}) and (\ref{eq_sg}).

As we identify $a_{-n}(\g)$ $(n >0, \g \in \G^*)$ with the $n$-th
power sum $p_n (\g)$ and the space $\SG$ with the space of
symmetric functions indexed by $ \G^*$, we may regard the Schur
function $s_{\lambda} (\g)$ associated to $\g$ and a partition
$\lambda$ a corresponding element in $\SG$. For $\lambda \in {\cal
P} (\G^*)$, we denote
\begin{eqnarray*}
  s_{\lambda} = \prod_{\g \in \G^*} s_{\lambda(\g)} (\g) \in \SG.
\end{eqnarray*}
Then $s_{\lambda}$ is the image under the characteristic map
$\ch$ of the character of an irreducible
representation $\chi^{\lambda}$ of $\Gn$ (cf. \cite{M}).

Denote by $c_n (c \in \G_*)$ the conjugacy class in $\Gn$ of
elements $(x, s) \in \Gn$ such that $s$ is an $n$-cycle and the
cycle product of $x$ lies in the conjugacy class $ c$. Denote by
$\sigma_n (c )$ the class function on $\Gn$ which takes value $n
\zeta_c$ (i.e. the order of the centralizer of an element in the
class $c_n$) on elements in the class $c_n$ and $0$ elsewhere. For
$\rho = \{ m_r (c) \}_{r \geq 1, c \in \G_*} \in {\cal P}_n
(\G_*)$, $\sigma_{\rho} = \prod_{r \geq 1, c \in \G_*} \sigma_r
(c)^{m_r (c)}$ is the class function on $\Gn$ which takes value
$Z_{\rho}$ on the conjugacy class of type $\rho$ and $0$
elsewhere. Given $\g \in R(\G)$, we denote by $\sigma_n (\g )$ the
class function on $\Gn $ which takes value $n \g (c) $ on elements
in the class $c_n, c \in \G_*$, and $0$ elsewhere.
%
%
%
%
%
We define an operator $ \widetilde{a}_{ -n} (\gamma), n >0$ to be
a map from $\RG$ to itself by the following composition
\[  R (\Gm) \stackrel{ \sigma_n ( \g ) \otimes }{\longrightarrow}
  R(\Gn) \bigotimes R (\Gm)  \stackrel{{Ind} }{\longrightarrow}
  R ( {\Gamma}_{n +m}).
\]
We also define another operator
$ \widetilde{a}_{ n} (\gamma), n >0$ to be a map from $\RG$
to itself (which is the adjoint of $\widetilde{a}_{ -n} (\gamma)$
with respect to the bilinear form (\ref{eq_series}) )
as the composition
\[
  R (\Gm)  \stackrel{ Res }{\longrightarrow}
   R(\Gn)\bigotimes R ( {\G }_{m -n})
 \stackrel{ \langle \sigma_n ( \g), \cdot \rangle }{\longrightarrow}
 R ( {\G }_{m -n}).
\]

The following theorem was established in \cite{W} (also see
\cite{FJW}).

\begin{theorem} \label{th_heis}
  The space $\RG$ affords a representation of the Heisenberg algebra
 $\hg $  by letting $ a_n (\g )$
 $( n \in \Bbb Z \backslash 0)$ act as $ \widetilde{a}_{ n} (\g )$,
 and $C$ as $1$. The characteristic map
 $\ch$ is an isomorphism of $\RG$
 and $\SG$ as representations over the Heisenberg algebra.
\end{theorem}
\section{Virasoro algebra and group convolution}  \label{sect_conv}
%
%
%
%
%
%
We first look at the case when $\G$ is trivial and so $\Gn$
becomes the symmetric group $S_n$. We simply write the $i$-th
power sum $p_i (\g)$ as $p_i$.

We consider the convolution product on the space of class
functions on $S_n$ with the class function $K_{(1^{n -2} 2)}$,
which takes value $1$ at elements of cycle type $ (1^{n -2} 2)$
and $0$ otherwise. It follows from (\ref{eq_charmap}) that $$ \ch
(K_{(1^{n -2} 2)} ) = \frac1{2 (n -2)!} p_1^{n -2} p_2 . $$ We
denote
\begin{eqnarray}   \label{eq_delta}
  \De = \hf \sum_{i, j \geq 1}
      \left( i j p_{i +j} \diffi \diffj
         + ( i +j) p_i p_j \diffij
      \right) .
\end{eqnarray}

Given $f, g \in R(S_n)$, it is natural to define the
convolution of two symmetric functions
$\ch (f)$ and $\ch (g)$ as
$$
  \ch (f) * \ch (g)  := \ch (f * g).
$$

The next theorem describes the effect of the convolution with
$K_{(1^{n -2} 2)}$ on the space of symmetric polynomials by means
of the characteristic map.
\begin{theorem} \label{th_symm}
 For any $f \in \C [S_n ]$, we have
 $$
   \ch (K_{(1^{n -2} 2)} * f) = \De \ch (f).
 $$
 Equivalently we have
 $$
   \frac 1{2 (n -2)!} p_1^{n -2} p_2 * \ch (f) = \De \ch (f).
 $$
\end{theorem}

\begin{demo}{proof}
 Take a transposition $(a,b)$ and a permutation $\tau$.

If $a$ and $b$ lie in different cycles of $\tau$, the product
$(a,b) \tau$ will have the effect of combining the two cycles, say
of cycle length $i$ and $j$ respectively, in $\tau$ containing
respectively $a$ and $b$ into a single cycle. For example, let
$n=7, a=3, b=5$ and $\tau = (1, 3) (2, 6, 5)(4, 7)$, then $(3, 5)
\tau = (1,5,2,6,3) (4, 7)$. Thus one $p_i p_j$ in $\ch (\tau )$ is
replaced by one $p_{i +j}$ in $\ch ( (a, b)\tau )$. Among all
transpositions  $(a, b) \in S_n$, there are exactly $i j$ of which
have the effect of replacing one $p_i p_j$ by one $p_{i +j}$.

On the other hand, if $a$ and $b$ lie in a same cycle, say of
cycle length $k$, of $\tau$, then the product $(a, b) \tau$ will
have the effect of splitting this cycle in $\tau$ containing $a$
and $b$ into two disjoint cycles. For example, let $n = 8,
a=3,b=5$ and $\tau = (8, 7, 3,1, 4, 5) (2, 6)$, then $(3, 5) \tau
= (8, 7, 5) (1,4, 3) (2, 6)$. More precisely the cycle of $\tau$
splits into two disjoint cycles of length $i$ and $k -i$ if $a$
and $ b$ in the cycle of $\tau$ are separated by $i-1$ and $k-i-1$
elements. Thus one $p_{k}$ in $\ch (\tau )$ is replaced by one
$p_i p_{k -i}$ in $\ch ( (a, b)\tau )$. We can see easily among
all possible transpositions there are $k$ (resp. $k/2$) of them
which have the effect of replacing one $p_{k}$ by one $p_i p_{k
-i}$ when $k/2 \neq i$ (resp. $k/2 = i$).

Combining the above considerations together, we have proved the
theorem.
\end{demo}

\begin{remark}  \rm
Our purpose of studying this convolution is to find a group
theoretic construction of the Virasoro algebra. It turns out that
such a convolution appeared earlier in a paper of Goulden
(\cite{Go}, Proposition~3.1) in his study of the number of ways of
writing permutations in a given conjugacy class as products of
transpositions.
\end{remark}

Recall that $s_{\la}$ denotes the Schur functions associated to
the partition $\la$. We denote by $f^{\la}$ the degree of
$s_{\la}$, i.e. the dimension of the irreducible representation of
$S_n$ corresponding to $s_{\lambda}$.

\begin{proposition}  \label{prop_reform}
  We have
\begin{eqnarray*}
\De (\schu ) = \frac {n (n -1)}{2 \f } \langle p_1^{n -2} p_2,
\schu \rangle \schu.
\end{eqnarray*}
\end{proposition}

\begin{demo}{Proof}
By the orthogonality relation of characters (\ref{eq_stand}) we
have $\langle s_\mu , s_\la \rangle = \delta_{\mu, \la}$, and thus
$p_1^{n -2} p_2 =\sum_{\mu} \langle p_1^{n -2} p_2, s_{\mu}
\rangle s_{\mu}.$ In addition Eq.~(\ref{eq_idem}) implies that
$s_{\mu} * s_{\la} = \frac{n!}{f^{\la}}\delta_{\mu,\la} s_{\la}.$
Thus by Theorem~\ref{th_symm}, we have
\begin{eqnarray*}
\De (\schu ) &=& \frac1{2 (n-2)!} p_1^{n -2} p_2 * \schu \\
 &=& \frac1{2 (n-2)!} \sum_{\mu} \langle p_1^{n -2} p_2, s_{\mu}
\rangle s_{\mu} * \schu \\
 &=&\frac {n (n -1)}{2 \f } \langle p_1^{n -2} p_2,
\schu \rangle \schu.
\end{eqnarray*}
\end{demo}
%
%
%
%
%
%
%

Now we return to the general case of wreath product $\Gn$. Given
$c \in \G_*$, we denote by $\laz \in {\cal P}_n (\G_*)$ the
function which maps $c$ to the one-part partition $(2)$, $c^0$ to
the partition $( 1^{n-2})$ and other conjugacy classes to $0$. We
denote by $K_{\laz}$ the sum of all the elements in the conjugacy
class corresponding to $\laz$.

We introduce the operator acting on $\RG$
$$
\De^{\g } =  \frac16 \int :\widetilde{a}^{\g} (z)^3 : z^2  dz
$$
where $ \widetilde{a}^{\g} (z)
 = \sum_{n \in \Z} \widetilde{a}_n (\g) z^{ -n -1}.$

We define
\begin{eqnarray} \label{eq_convop}
 \Dec = \sum_{\beta \in \G^*}
  \frac{|\G|^2 \beta (c^{-1})}{d_{\beta}^2 \zeta_c}
  \Delta^{\beta}.
\end{eqnarray}
The operator $\Dec$ reduces to (\ref{eq_delta}) for $\G$
trivial. By the characteristic map, the operator $\De^{\g}$ can be
identified with the differential operator $$  \hf \sum_{i, j \geq
1}
      \left( i j p_{i +j}(\g) \diffig \diffjg
         + ( i +j) p_i(\g) p_j(\g) \diffijg
      \right).
$$
Note that
$$
 \ch (K_{\laz}) =
  \frac{1}{Z_{\laz }} p_1 (c^0)^{n -2} p_2 (c),
$$
where
\begin{eqnarray}  \label{eq_order}
 Z_{\laz} = 2 (n -2)! | \G |^{n -2} \zeta_c
\end{eqnarray}
is the order of the centralizer of an element in the conjugacy
class associated to ${\laz}$, cf. (\ref{eq_centr}). The following
theorem generalizes Theorem~\ref{th_symm}.

\begin{theorem} \label{th_main}
 Given $f \in  \RG$, we have
 $$
  \ch (K_{\laz} * f) = \Dec \ch (f) .
 $$
 Equivalently we have
 $$
  \frac{1}{Z_{\laz }} p_2 (c) p_1(c^0)^{n -2} * \ch (f)
  = \Dec \ch (f).
 $$
\end{theorem}

We need some preparation for the proof of the theorem. Recall that
the image of the irreducible character associated to $\la \in
{\cal P}_n (\G^*)$ under the characteristic map is $\schu =
\prod_{\g} \sge$. We set $n_{\g} = |\la (\g)|$. The first lemma
below is straightforward.

\begin{lemma}  \label{lem_zero}
 Given $\be \in \G^*$ and a sequence of non-negative
 integers $m_{\g}, \g \in \G^*$ such that
 $m_{\g} \neq n_{\g}$ for at least one $\g $,
 we have
 \begin{eqnarray*}
   \left\langle p_1 (\be)^{m_{\be} -2} p_2 (\be)
        \prod_{\g \neq \be} p_1 (\g)^{m_{\g}} ,
        \schu
   \right\rangle  = 0 .
 \end{eqnarray*}
\end{lemma}

\begin{lemma}  \label{lem_comp}
Given $\be \in \G^*$, $\la \in {\cal P}_n (\G^*)$, and setting
$n_{\g} = |\la (\g)|$ (for all $\g \in \G^*$), we have
 \begin{eqnarray*}
  & & \frac {n_{\be} (n_{\be} -1)}{2}
  \left\langle p_1 (\be)^{n_{\be} -2} p_2 (\be)
        \prod_{\g \neq \be} p_1 (\g)^{n_{\g}} ,
        \schu
  \right\rangle \schu      \\
  & =&\left( \prod_{\g} \fge
       \right)
       \De^{\be} \sbe \cdot
       \prod_{\g \neq \be} \sge.
 \end{eqnarray*}
\end{lemma}

\begin{demo}{Proof}
 It follows from
 Proposition~\ref{prop_reform} that
 \begin{eqnarray}  \label{eq_factor}
  \De^{\be} (\sbe ) =
  \frac {n_{\be} (n_{\be} -1)}{2 \fbe}
  \left\langle p_1 (\be)^{n_{\be} -2} p_2 (\be), \sbe
  \right\rangle \sbe.
 \end{eqnarray}

    Noting also that
 \begin{eqnarray}   \label{eq_deg}
   \left\langle p_1 (\g)^{n_{\g}}, \sge  \right\rangle = \fge,
 \end{eqnarray}
 we calculate that
 \begin{eqnarray*}
  & & \frac {n_{\be} (n_{\be} -1)}{2}
  \left\langle p_1 (\be)^{n_{\be} -2} p_2 (\be)
        \prod_{\g \neq \be} p_1 (\g)^{n_{\g}} ,
        \schu
  \right\rangle \schu                       \nonumber   \\
  & =& \frac {n_{\be} (n_{\be} -1)}{2}
   \left\langle p_1 (\be)^{n_{\be} -2} p_2 (\be), \sbe
   \right\rangle \sbe \cdot                        \nonumber  \\
   & & \cdot \prod_{\g \neq \be}
       \left\langle p_1 (\g)^{n_{\g}}, \sge
       \right\rangle \sge                       \nonumber   \\
  & =& \fbe  \De^{\be} \sbe \cdot
       \prod_{\g \neq \be} \fge \sge
        \quad \mbox{ by Eqs. } (\ref{eq_factor})
             \mbox{ and }(\ref{eq_deg}),
                                 \nonumber   \\
  & =& \left( \prod_{\g \in \G^*} \fge
       \right)
       \De^{\be} \sbe \cdot
       \prod_{\g \neq \be} \sge.
 \end{eqnarray*}
 This proves the lemma.
\end{demo}

\begin{demo}{Proof of Theorem~\ref{th_main} }
Consider the irreducible character associated to $\la \in {\cal
P}_n (\G^*)$ which maps by the characteristic map $\ch$ to $\schu
= \prod_{\g} \sge$. We set $n_{\g} = |\la (\g)|$.

 We calculate that
 \begin{eqnarray}   \label{eq_expan}
  & & \frac 1{2 (n -2)!} p_1(c^0)^{n -2} p_2 (c)   \nonumber \\
  & =& \frac 1{2 (n -2)!}
    \left( \sum_{\g \in \G^*} d_{\g} p_1 (\g)
    \right)^{n -2}
    \left(\sum_{\be \in \G^*} {\be}(c^{-1}) p_2 (\be)
    \right)  \quad \mbox{ by Eq.}~(\ref{eq_change}), \nonumber  \\
  & =& \hf
   \left(
       \sum_{ \{ m_{\g}\} }
       \prod_{\g  \in \G^*}
       \frac{ (d_{\g} p_1 (\g))^{m_{\g}} }{m_{\g}!}
   \right)
    \left(\sum_{\be  \in \G^*} {\be}(c^{-1}) p_2 (\be)
    \right)    \nonumber \\
  & =& \hf \sum_{ \{ m_{\g}\} } \sum_{\be \in \G^*}
    {\be}(c^{-1}) p_2 (\be)
    \frac{(d_{\be} p_1 (\be))^{m_{\be} }}{(m_{\be})!}
    \prod_{\g \neq \be}
       \frac{ (d_{\g} p_1 (\g))^{m_{\g}} }{m_{\g}!} ,
 \end{eqnarray}
 where the sum $\sum_{ \{ m_{\g}\} }$ is taken over all
 sequences of nonnegative integers $m_{\g}, \g \in \G^*$
 such that $\sum_{\g} m_{\g} =n -2$.

 Recall \cite{M, Z} that the degree of the irreducible character
 associated to $\lambda \in {\cal P}_n (\G^*)$ is
 \begin{eqnarray}   \label{eq_total}
  \f = n ! \prod_{\g} d_{\g}^{n_{\g}} \fge / n_{\g} !.
 \end{eqnarray}
 We also obtain by using (\ref{eq_idem})
 and $|\Gn| = n! |\G|^n$ that

 \begin{eqnarray}  \label{eq_apply}
  && \frac{\f}{n! | \G |^n}
   \left(
       p_2(\be) p_1(\be)^{n_{\be} -2}
       \prod_{\g \neq \be} p_1(\g)^{n_{\g}}
         \right)  * \schu       \nonumber   \\
   & =&
   \left \langle
       p_2(\be) p_1(\be)^{n_{\be} -2}
       \prod_{\g \neq \be} p_1(\g)^{n_{\g}} ,  \schu
   \right \rangle \schu.
 \end{eqnarray}

 We then have
 \begin{eqnarray*}
  &  & \frac 1{Z_{\laz}}
        p_1(c^0)^{n -2} p_2 (c) * \schu  \\
  & =& \frac{1}{2 |\G|^{n -2} \zeta_c}\sum_{\beta \in \G^*}
    {\be}(c^{-1}) p_2 (\be)
    \frac{(d_{\be} p_1 (\be))^{n_{\be} -2}}{(n_{\be} -2)!}
    \prod_{\g \neq \be}
       \frac{ (d_{\g} p_1 (\g))^{n_{\g}} }{n_{\g}!}  * \schu  \\
  & &   \quad \quad \quad \quad
       \mbox{ by Eqs.}~(\ref{eq_order}), (\ref{eq_expan}),
       \mbox{ and Lemma}~\ref{lem_zero}, \\
  & =& \frac{1}{ |\G|^{n -2}\zeta_c }
      \prod_{\g} \frac{d_{\g}^{n_{\g}}}{n_{\g}!}
      \sum_{\be \in \G^*}
      \frac{n_{\be} (n_{\be} -1) \be(c^{-1})}{2 d_{\be}^2}
     \left(
       p_2(\be) p_1(\be)^{n_{\be} -2}
       \prod_{\g \neq \be} p_1(\g)^{n_{\g}}
         \right)  * \schu  \\
  & =& \frac{\f}{n!|\G|^{n -2}\zeta_c}  \sum_{\be \in \G^*}
       \frac{n_{\be} (n_{\be} -1) \be(c^{-1})}{2 d_{\be}^2
        (\prod_{\g} \fge)}
     \left(
       p_2(\be) p_1(\be)^{n_{\be} -2}
       \prod_{\g \neq \be} p_1(\g)^{n_{\g}}
         \right)  * \schu \\
         &&\qquad \qquad \qquad \qquad \mbox{ by Eq.}~(\ref{eq_total}), \\
  & =& \frac{|\G|^2}{\zeta_c}  \sum_{\be \in \G^*}
       \frac{n_{\be} (n_{\be} -1) \be (c^{-1})}{2 d_{\be}^2(\prod_{\g} \fge)}
     \left\langle
       p_2(\be) p_1(\be)^{n_{\be} -2}
       \prod_{\g \neq \be} p_1(\g)^{n_{\g}} ,
         \schu \right\rangle s_{\la} \\
         &&\qquad \qquad \qquad \qquad  \mbox{ by Eq.}~(\ref{eq_apply}), \\
  & =&  \sum_{\be \in \G^*}
       \frac{|\G|^2\be (c^{-1})}{d_{\be}^2 \zeta_c}
     \left(
             \De^{\be} \sbe \cdot
       \prod_{\g \neq \be} \sge
         \right)  \\
         &&\qquad \qquad \qquad \qquad  \mbox{ by Lemma}~\ref{lem_comp}, \\
  & =& \Dec \schu.
\end{eqnarray*}
The last equation holds since $\De^{\be} \sge = 0$ for $\g \neq
\be$. This finishes the proof, since $\schu, \lambda \in {\cal P}
(\G^*)$ form a basis of the space $\RG$.
\end{demo}

\begin{corollary}
For the identity conjugacy class $c^0 \in \G_*$, we have
$\De_{c^0} = \sum_{\be \in \G^*} \frac{|\G|}{d_{\be}} \De^{\be}.$
\end{corollary}

Combining (\ref{eq_convop}), Lemma~\ref{lem_ham},
Theorem~\ref{th_heis} and Theorem~\ref{th_main}, we obtain the
following.

\begin{theorem}
The operator $\Dec$ acting on $\RG$ is realized by the convolution
with $K_{\laz}$. The commutation relation between $\Dec$ and the
Heisenberg algebra generator $\widetilde{a}_{n} (\g)$ (constructed
in a group theoretic manner) is given by
 \[
 [ \Dec, \widetilde{a}_{n} (\g)]
  = -\frac{n|\G |^2 \g (c^{-1})}{d_{\g}^2 \zeta_c} L_{n}^{\g},
\]
where the operators $L_n^{\g}$ acting on the space $\RG$ satisfy
the Virasoro commutation relations:
 \[
    [ L_n^{\be}, L_m^{\g} ]
    = \delta_{\be  \g} L_{n +m}^{\g} +
     \frac{n^3 -n}{12} \delta_{\be  \g} \delta_{n, -m}.
 \]
\end{theorem}

Frenkel: Department of Mathematics, Yale University, New Haven, CT 06520;
Wang: Department of Mathematics, North Carolina State University,
Raleigh, NC~27695, wqwang@math.ncsu.edu.

\end{document}